\documentclass[11pt]{article}
\usepackage{float}
\usepackage{url}

\usepackage{fullpage, amsmath, amssymb, amsthm, epsfig, color, graphicx}
\newenvironment{@abssec}[1]{%
     \if@twocolumn
       \section*{#1}%
     \else
       \vspace{.05in}\footnotesize
       \parindent .2in
         {\bfseries #1. }\ignorespaces
     \fi}
     {\if@twocolumn\else\par\vspace{.1in}\fi}
\newenvironment{keywords}{\begin{@abssec}{Key words}}{\end{@abssec}}
\newenvironment{AMS}{\begin{@abssec}{AMS subject classification}}{\end{@abssec}}

\def\eqbd{\mathop{{:}{=}}}

\def\Re{\mathop{\rm Re}\nolimits}
\def\openC{{\rm C\kern-.48em\vrule width.06em height.6em depth-.02em
                 \kern.48em}}
\def\openQ{{{\rm Q\kern-.21cm\vrule width.6pt height 6.2ptdepth-.2pt \kern.21cm}}}
\def\openR{{{\rm I}\kern-.16em {\rm R}}}
\def\openZ{{{\rm Z}\kern-.28em{\rm Z}}}
\def\openT{{{\rm T}\kern-.42em {\rm T}}}
\def\openH{{{\rm I}\kern-.16em {\rm H}}}
\def\openK{{{\rm I}\kern-.16em {\rm K}}}
\def\openL{{{\rm I}\kern-.16em {\rm L}}}
\def\openM{{{\rm I}\kern-.16em {\rm M}}}
\def\openN{{{\rm I}\kern-.16em {\rm N}}}
\def\openP{{{\rm I}\kern-.16em {\rm P}}}

\def\eqbd{\mathop{{:}{=}}}

\let\N\openN

\def\Z{\mathbb Z}

\def\belowrightarrow#1{{{{}\over\ #1\ }\kern-1.1em\to}}
\def\l2{{L_2}}

\def\nt{\noindent}

\def\eqbd{\mathop{{:}{=}}}

\def\Re{\mathop{\rm Re}\nolimits}

\def\eqbd{\mathop{:}{=}}

\def\belowrightarrow#1{{{{}\over\ #1\ }\kern-1.1em\to}}
\def\l2{{L_2}}

\begin{document}

\title{Classical values of Zeta, as simple as possible but not simpler}
\author{ Olga Holtz \\[1mm] {\small Department of Mathematics, University of California, Berkeley }}
\date{\small March 24, 2024}
\maketitle

\nt

\nt

\begin{abstract} \nt
This short note for non-experts means to demystify the tasks of evaluating the Riemann Zeta Function at nonpositive integers and at even natural numbers, both initially performed by Leonhard Euler.  Treading in the footsteps of G. H. Hardy and others, I re-examine Euler's work on the functional equation for the Zeta function, and explain how both the functional equation and all `classical' integer values can be obtained in one sweep using only Euler's favorite method of generating functions. As a counter-point, I also present an even simpler argument essentially due to Bernhard Riemann, which however requires Cauchy's residue theorem, a result not yet available to Euler. As a final point, I endeavor to clarify how these two methods are organically linked and can be taught as an intuitive gateway into the world of Zeta functionology.
\end{abstract}

\section{History in a nutshell}

The notorious Basel problem, posed by Pietro Mengoli in 1650, asked for the precise evaluation of the infinite  sum
$$\sum_{n=1}^\infty {1\over n^2}={1\over 1^2}+{1\over 2^2}+{1\over 3^2}+\cdots ,$$  
which a modern reader would instantly recognize as the value of the celebrated Riemann Zeta function at the point $2$. This problem was solved by Leonhard Euler by 1734, with his solution read at The Saint Petersburg Academy of Sciences on December 5, 1734.

Euler's investigation of this and related problems began  in 1731 and continued at least till 1749. Euler's {\em tour de force\/} led to his solution or, in fact, multiple solutions to the Basel problem, followed by his successful evaluation of $\zeta$ at all positive even integers and at all negative integers, his infinite product representation for $\zeta$, and his lesser-known formulation of its functional equation, albeit only with a partial proof. The fascinating history of Euler's work on this subject is laid out in~\cite{Ayoub1974}.

About a century later in his groundbreaking memoir~\cite{riemann1859uber} (for an English translation, see~\cite[pp.~190--198]{borwein2008}) Bernhard Riemann re-introduced $\zeta$ as a function of a complex variable $s$, established a number of its integral representations,  found and justified its analytic continuation beyond the half-plane $\Re s >1$, proved its functional equation, and formulated The Riemann Hypothesis. Riemann's memoir must have been the most influential nine pages of mathematics ever written.

\section{Teaching orthodoxy and heterodoxy}

These days, undergraduate students are introduced to the values of the $\zeta$-function at positive even integers and all negative integers (which I will call `classical values' here) roughly as follows. First, the value $\zeta(2)$ is derived, usually using a partial fraction decomposition of the cotangent function, which, in turn, is justified via the so-called Herglotz trick~\cite[Chapter~26]{Aigner2018} or by using residues if the students are familiar with complex analysis. The same method also yields the values of $\zeta$ at all positive even integers. 

The $\zeta$-function of a complex variable $s$  is then defined for $\Re s>1$ using one of its integral forms and continued outside that half-plane. This continuation is usually performed using the Jacobi theta function and its functional equation, which, in turn, requires a proof of  Poisson summation, a mysterious formula with a demanding proof for undergraduates. Once a symmetric integral form of $\zeta (s)$ is proved, its functional equation follows. Then the values of $\zeta$ at negative integers are obtained using the functional equation. If the students are not familiar with complex analysis, the key facts of this derivation are merely stated.

This derivation of all `classical' values of $\zeta$ is lengthy, and, to my mind, not particularly edifying. 

The goal of this note is to present two simple ways to obtain all classical values of~$\zeta$. Both arguments will proceed in the order I wish I had learned as an undergraduate myself: from evaluating $\zeta$ at all negative integers, to its functional equation, to its evaluation at all even positive integers. I will argue below that this `counterintuitive' order is in fact very `natural'. I will also endeavour to strip each of its steps from all complexity.

The first way is described in Section~\ref{sec:riemann} following Riemann's approach in his memoir~\cite{riemann1859uber}. This approach requires only basic complex analysis, no Poisson summation formula, and, in fact, no Jacobi theta function.  The second way presented in Section~\ref{sec:euler} streamlines Euler's approach in his series of works~\cite{euler1731,euler1732/33, euler1734, euler1736,euler1737,euler1740,euler1749}. This approach does not even require complex analysis, only results and techniques available to Euler in 1730--1750.

Interestingly, apart from exactly one sentence about the values of $\zeta$ being zero at negative even integers, Riemann himself did not even address the classical values of $\zeta$ in his memoir -- or in his Nachlass, according to Siegel~\cite{siegel1931, barkan2018}. Nevertheless, I have no doubt that the approach spelled out in Section~\ref{sec:riemann} below was clear to Riemann. I am equally certain that the simplified version of Euler's arguments presented in Section~\ref{sec:euler} would have been clear to Euler. 

Thus, I claim no results below as new. However, their arrangement must be new. 

The first impetus for this note was my dismay at my own arduous answer to the question ``How to evaluate $\zeta$ at all classical points?'' once asked by a curious undergraduate. The second impetus was my disturbing realization, while answering this question, that I did not {\em really\/} know (!) why the Bernoulli numbers arise at the integer values of $\zeta$. This is my attempt to redeem myself. 

One more thing before we begin. The  reader might wonder how much Riemann actually knew about Euler's work on `his' $\zeta$-function. Riemann's expert knowledge of Latin and French (among other languages) and his remarkable ability to absorb mathematical knowledge by quick reading are well documented~\cite{riemannbio}. This, together with an assurance of Andr\'e Weil given to Raymond Ayoub as reported in~\cite{Ayoub1974}, makes me very certain that Riemann was intimately familiar with Euler's research on the $\zeta$-function, most likely including Euler's results on its functional equation.  I should add that Andr\'e Weil~also made an intriguing conjecture about a possible influence of Eisenstein's unpublished proof of the functional equation for the Dirichlet $L$-function modulo $4$ on Riemann's epic paper~\cite{riemann1859uber}, along with a plausible explanation how Eisenstein's copy of Gauss's {\em Disquisitiones\/} with that proof written as an annotation reached the mathematical library in Giessen (see~\cite{Weil1987, Weil1989}).


\section{Riemann's Way}\label{sec:riemann}

\subsection{Heart of the matter}

Riemann's starting point is the following observation based on a change of variables:
$$ \int_0^\infty e^{-nx} x^{s-1} dx = n^{-s} \int_0^\infty e^{-x} x^{s-1} dx = n^{-s} \Gamma (s),$$
with the last equality arising from the integral definition of the Gamma function due to Euler.

Riemann sums such integrals from $1$ to $\infty$, obtaining
$$\Gamma(s)\sum_{n=1}^\infty  n^{-s}= \Gamma(s) \zeta(s) =\int_0^\infty x^{s-1} (e^{-x} +e^{-2x}+\cdots) \, dx = 
\int_0^\infty {x^{s-1} e^{-x} \over 1- e^{-x}} \, dx = \int_0^\infty  { x^{s-1} dx  \over e^{x} -1}.$$
In a striking turn, Riemann then switches the sign of $x$ in the numerator to $-x$ and the contour of integration from the half line $(0,+\infty)$ (or $[0,+\infty)$) to the so-called Hankel contour, which I will denote by $\gamma$. Incidentally, Riemann could have flipped that contour around the origin instead, sending it to infinity leftwards, but he proceeds as above, obtaining
\begin{equation}  \int_\gamma {(-x)^{s-1} dx \over e^x -1}. \label{zetacontour} \end{equation}

\begin{figure}[h]
  \centering
  \includegraphics[width=2.3in]{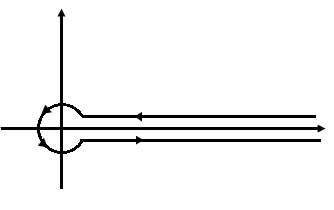}
  \caption{The Hankel contour $\gamma$ (courtesy of S. Westerway via Wikipedia).}
  \label{fig:Hankel_Pic}
\end{figure}
\nt
This integration presupposes that the (usually) multi-valued function $(-x)^{s-1}=e^{(s-1)\log (-x)}$
is classically defined, so that the logarithm of $-x$ is real when $x$ is negative. The point of this clever trick is to `extract' the  value of $\zeta(s)$ for any complex value $s$ from the integral~(\ref{zetacontour}):
\begin{equation}  \int_\gamma {(-x)^{s-1} dx \over e^x -1} =  (e^{-\pi s i} -e^{\pi s i})  \int_0^\infty {x^{s-1} dx \over e^x -1}
=-2i \sin (\pi s)\, \Gamma(s) \zeta(s).  \label{zetagamma} \end{equation}
In fact, just these few lines of Riemann's memoir provide a key to the $\zeta$-values at all nonpositive integers. Indeed, if $s$ is a nonpositive integer, the function $(-x)^{s-1}$ is single-valued and meromorphic with the only pole at the origin. So the contour $\gamma$ can be deformed into a loop around the origin without changing the value of the integral~(\ref{zetacontour}). That value can be now obtained, by Cauchy's residue theorem, from the coefficient of $x^{-1}$ in the Laurent series of the integrand, i.e., from the coefficient of $x^{-s+1}$ in the Taylor-Maclaurin series for $x/(e^x-1)$. (Remember that $-s$ is a nonnegative integer.) And the latter series was already computed by Euler (of course!),
who showed that
$${x\over e^x-1} = \sum_{n=0}^\infty {B_n \over n!}  x^n,$$
where $(B_n)$ denotes  the sequence of Bernoulli numbers. 

A remaining technicality is to make sense of the value $\sin (-\pi n) \Gamma(-n)$, which appears in the right-hand side of~(\ref{zetagamma}).
Since the Gamma function has poles at all nonpositive integers, this is an indeterminate of the sort $0\cdot \infty$, but it can be resolved simply, using the defining multiplicative property of the Gamma function $n$ times:
$$\lim_{x\to -n} \sin (\pi x ) \,\Gamma(x) =\lim_{x\to -n} {\sin (\pi x) \, \Gamma(x+n+1) \over (x+n) (x+n-1)\cdots (x+1)x} =
{(-1)^n \pi \, \Gamma(1) \over (-1) (-2)\cdots (-n)} = {\pi \over n!}.$$
Recalling that the integrand of~(\ref{zetacontour}) contains the sign $(-1)^{n-1}$, we now get
$$ 2\pi i \; {(-1)^{n-1}B_{n+1} \over (n+1)!}= -2i \sin (-\pi n) \Gamma(-n) \zeta(-n) =  - 2 i \, {\pi\over n!} \, \zeta(-n),$$
and hence \begin{equation} \zeta(-n)={(-1)^n B_{n+1} \over n+1}. \label{zetanegative} \end{equation}
This formula holds for all nonpositive integers $-n$. 

 Riemann in~\cite{riemann1859uber} immediately proceeds to note that the integration in~(\ref{zetacontour}) can be turned `inside out' whenever $\Re s<0$. This amounts to augmenting the Hankel contour $\gamma$ with a big circle $|z|=R$, whose  contribution to the integral will be zero when $R\to \infty$. 
 Riemann does not elaborate on this point, which can be seen by multiplying the arc length $2\pi R$ by a bound of order $O(R^{\Re s -1})$ on the integrand, which we get if the contour stays away from the poles on the imaginary axis, say, by taking $R$ to be an odd multiple of $\pi$.  The resulting bound, of order $O(R^{\Re s})$, tends to $0$.
 
 In that way, the contour~$\gamma$ can be thought of carving a slice out of the complex plane, and Cauchy's residue theorem can be used not at the origin but at all remaining poles of the integrand, i.e., the values $2\pi n i$, $n\in \Z\setminus \{0\}$.  The resulting contribution from each pole is then $(-2\pi i) (-2\pi n i)^{s-1} $, so the total is
$$ 2 \sin (\pi s) \Gamma (s) \zeta(s) = (2\pi)^s \sum_{n=1}^\infty n^{s-1} ((-i)^{s-1} +i^{s-1}) =  (2 \pi)^s  \zeta(1-s)  \, 2\sin (\pi s/2) ,$$
which simplies to
\begin{equation} 2 \cos (\pi s/2) \Gamma(s) \zeta(s)  = (2\pi)^s  \zeta(1-s). \label{funceq}  \end{equation}
This is one of the simplest forms of the famous functional equation for the $\zeta$-function.

Now apply~(\ref{funceq}) to evaluate $\zeta$ at positive even integers. We have
$$     2 \cos (\pi (1-2n)/2) \Gamma (1-2n) \zeta (1-2n) = (2\pi)^{1-2n} \zeta(2n)  $$
and so must only make sense of the term $ 2 \cos (\pi(1-2n)/2) \Gamma(1-2n)$, which we already know how to do since the limiting value of that 
product is equal to the limiting value of the product $ \sin (1-2n) \Gamma(1-2n)$ modulo the sign $\sin ((1-2n) \pi/2)=(-1)^n$. Thus, we get
$$\zeta(2n)= (2\pi)^{2n-1} {(-1)^n  \pi \over (2n-1)!} \,\zeta(1-2n)= {(2\pi)^{2n-1} (-1)^n \pi \over (2n-1)!} \cdot{(-1)^{2n-1} B_{2n} \over 2n} =
 { (-1)^{n-1} (2\pi)^{2n} B_{2n} \over 2 (2n)!}.$$
 That's it! But the reflection formula won't yield the values of $\zeta$ at positive odd integers. [Why?]

\section{Euler's Way}\label{sec:euler}

Next, let us imagine how Euler might have implemented the same program: 1) evaluate $\zeta$ at the nonnegative integers, 2) establish its functional equation, and 3) use it to evaluate $\zeta$ at the even positive integers. While indulging this fantasy, let us agree to use only results and techniques known to Euler around the time he worked on the Basel problem and its generalizations.

\subsection{Heart of the matter}\label{sec:euler_heart}

Euler was a virtuoso of generating functions. So let us start with the generating function for all nonpositive integer values of $\zeta$.  If that task appears daunting, gentle reader, fear not!

Let $S_m(n)$ denote the finite sum $\sum_{k=1}^n k^m$, and consider the exponential generating function $G(z,n)$ for the sequence $(S_m(n))_{m=0}^\infty$. That means
$$G(z,n)\eqbd \sum_{m=0}^\infty S_m(n) {z^m \over m!} = \sum_{m=0}^\infty \sum_{k=1}^n { (kz)^m \over m!} = \sum_{k=1}^n e^{kz} =
e^z \cdot {1-e^{nz}\over 1-e^z}={1-e^{nz} \over e^{-z} -1}.$$
If $n$ tends to $\infty$, the sum $S_m(n)$ must be replaced by the value of the $\zeta$-function at $-m$ and the generating function by its limit $1/(e^{-z}-1)$. An `unpleasant' term $-1/z$ appears in the Laurent expansion of this function at the origin, which we will, for now, unceremoniously discard. We get
$$ \sum_{m=0}^\infty \zeta(-m) {z^m\over m!} = {1\over e^{-z}-1} +{1\over z} =  \sum_{m=0}^\infty (-1)^{m-1}B_m {z^{m-1} \over m!} +{1\over z}.$$
Comparing coefficients, we get: $\zeta(-m) =(-1)^m B_{m+1}/(m+1)$, which is formula~(\ref{zetanegative}).

Now consider the odd version of the above generating function:
$${e^{-z} +1 \over e^{-z} -1}={1\over e^{-z}-1}+{ e^{-z} \over e^{-z} -1}={1\over e^{-z}-1}-{1\over e^z -1} =-{2\over z} +
2\sum_{m=0}^\infty \zeta(-2m-1) {z^{2m+1} \over (2m+1)!}  .$$ On the  other hand,
$${e^{-z} +1 \over e^{-z} -1} = {e^{-z/2}+e^{z/2} \over e^{-z/2}-e^{z/2}} =  {2\cos (-iz/2) \over -2i \sin(-iz/2)}= i \cot (-iz/2),$$
and the latter function has another, magical, representation, due to Euler (who else?)~\cite{Euler1748}:
\begin{equation} \pi \cot (\pi x) = {1\over x} + \sum_{n=1}^\infty  \left( {1\over x+n } +{1\over x-n}\right). \label{cotan} \end{equation}
Replacing $\pi x$ by  $-iz/2$ in~(\ref{cotan}), we now get
\begin{eqnarray*} {e^{-z} +1 \over e^{-z} -1}  =   {i\over \pi} \left[ {2\pi i \over z} +  \sum_{n=1}^\infty \left( 
{1\over {z\over 2\pi i} +n}  + {1 \over {z\over 2\pi i} -n}    \right)  \right] 
= -{2\over z}-2 \sum_{n=1}^\infty \left( {1\over z+2\pi i n} +  {1\over z-2\pi i n} \right)
  \end{eqnarray*}
  Therefore,
\begin{eqnarray*} 
 {e^{-z} +1 \over e^{-z} -1} +{2\over z} &=&  -\sum_{n=1}^\infty {4 z \over z^2 + 4\pi^2 n^2}  
 =   - \sum_{n=1}^\infty {z \over \pi^2 n^2} \cdot {1\over 1+\left( {z\over 2\pi n} \right)^2} =\sum_{m=0}^\infty {z^{2m+1} \over \pi^2}
 \sum_{n=1}^\infty {(-1)^{m+1} \over n^{2m+2} (2\pi)^{2m}} .
\end{eqnarray*}
This yields a  simple version of the functional equation~(\ref{funceq}), which is enough for our purposes:
$$ {2\zeta(-2m-1)\over (2m+1)!} = {(-1)^{m+1} \zeta(2m+2)  \over 2^{2m} \pi^{2m+2}},$$
which finally implies $\zeta(2m)={ (-1)^{m-1} (2\pi)^{2m} B_{2m} \over 2 (2m)!}$ for all natural numbers $m$, Q.E.D.

\subsection{More history}\label{sec:history}

A few comments are in order to defend my contention that the preceding argument could have been obtained by Euler some time between 1730 and 1750. At that time, Euler was well familiar with the generating function for the Bernoulli sequence and could certainly perform all generating function manipulations shown above up to and including  the cotangent formula~(\ref{cotan}). 

The timing of the latter formula is a bit tricky, since it was officially published by Euler only in 1748 in~\cite{Euler1748}, which is however still within the time period considered. On the other hand, the cotangent function is simply the logarithmic derivative of the sine function, and Euler mentioned the infinite product for the sine function already in his first solution to the Basel problem in 1734, so it would not be a stretch to say that the formula~(\ref{cotan}) was known to Euler already in 1734.

As to the functional equation for the $\zeta$-function, Euler arrived at it by 1749. By that time, he long knew the positive even integer values of $\zeta$ and also understood how to evaluate $\zeta$ at negative integers. By comparing the values of $\zeta$ at the positive and negative integers, he obtained, {\it a posteriori}, the functional equation for all integer values of $\zeta$, which he  then naturally conjectured to hold true for all real values. Euler was also able to verify its validity by taking the limit at~$1$ (where $\zeta$ has a pole), and (numerically) check it at the half-integer points. Euler's partial proof of the functional equation for the $\zeta$-function  was published  in \cite{euler1749} in 1749, 110 years before Riemann's memoir~\cite{riemann1859uber} on the subject. Remarkably, both Euler's and Riemann's memoirs were submitted to the Berlin Academy of Sciences!

Detailed discussions of Euler's 1749 paper were written by G.H. Hardy in \cite{Hardy1963} and R. Ayoub in~\cite{Ayoub1974}, with a shorter discussion by A. Weil in \cite{Weil1987}. Peculiarly, Hardy does not actually comment on the way Euler obtained the negative integer values of $\zeta$ in \cite{euler1749} but Ayoub does. 

A key observation that led Euler to his understanding of the values of $\zeta$ at negative integers is the famous Euler-MacLaurin formula. It was initially claimed by Euler already in 1732/33 \cite{euler1732/33} and proved by him in 1736 in \cite{euler1736}. 
His derivation of it is strikingly modern.  Let's look at it right away. 

Given any function $f$ and a constant $a$, Euler defines its `summatory' function
$$ S(x) \eqbd \sum_{n=0}^\infty f(x+na),$$  then rewrites this relation as 
$ -f(x) = S(x+a)-S(x)$
and expands $S(x+a)$ into a Taylor series at $x$ to obtain:
$$ -f(x) = \sum_{n=1}^\infty {a^n S^{(n)}(x) \over n!}.$$ 
Using the modern symbol $D$ for the differentiation operator $d/dx$, Euler (and we) can rewrite it as
$$ -f(x) = \left(I +{aD\over 2!} + {a^2 D^2 \over 3!} + \cdots\right) a S^{(1)}(x) = \left(  e^{aD} - I \over D  \right) S^{(1)}(x),$$
whose invertion yields
$$ S^{(1)} (x) = - \left( { D \over  e^{aD} - I }\right)f(x),$$
which is exactly the generating function of the Bernoulli sequence (in the operator $aD$)! This produces the  Euler-MacLaurin formula:
\begin{equation} S^{(1)}(x) = - \sum_{n=0}^\infty  {a^{n-1} B_n \over n!} f^{(n)}(x). \label{EM} \end{equation}

\subsection{Less mystery}\label{sec:mystery}

But, but, but... you must be asking yourself, isn't the `derivation' in Section~\ref{sec:euler_heart} not quite kosher? In fact,  wildly non-rigorous? Let's take a closer look.

The first problematic moment which the perceptive reader no doubt noticed right away is the transition between the generating function 
$G(z,n)$ that encodes all finite sums $S_m(n)=\sum_{k=1}^n k^m$ and the generating function $G(z)$ that is supposed to encode the infinite sums
$\zeta(-m)$, $m=0,1,\ldots$. And, before we can address this transition, we are faced with the fact that each sum $S_m(n)$
in fact diverges as $n\to \infty$! So in what sense can we possibly think of each sum $S_m(n)$ as approaching $\zeta(-m)$? And how, for that matter, did Euler  approach those  values of zeta?

After all, Euler's  contemporaries were wary or downright scared of divergent series. Even much later, at the time of Hardy's writing his book Divergent Series \cite{Hardy1963}, that title alone  would scandalize other mathematicians, so Hardy could be thought of as `trolling' his colleagues. However, Euler was not easily disturbed by divergence. Many of Euler's derivations used divergent series with great flair, arriving at correct results.  So in what sense did Euler think of `summing' quantities 
like $S_m(n)$ for $n\to \infty$? And what does it have to do with the Euler-MacLaurin summation?

Unsurprisingly, Euler found a brilliant if slightly roundabout way to access the values of $\zeta$ at nonpositive integers. Of course, he was perfectly aware that each sum $S_m(n)$ tends to $+\infty$ when $n$ does. Euler's inventive solution was to switch to their  alternating counterparts 
$\sum_{k=1}^n (-1)^{k+1} k^m$ instead.  The skeptical reader might now ask: why? or rather, what for?

To understand this, let us indulge in a bit of wishful thinking. Imagine that both series $\sum_{k=1}^\infty k^m$ and $\sum_{k=1}^\infty (-1)^{k+1} k^m$ were convergent. Then we would have
$$\sum_{k=1}^\infty  (2k)^m = 2^m \sum_{k=1}^\infty k^m,$$ therefore
\begin{equation} \sum_{k=1}^\infty (-1)^{k+1} k^m = \sum_{k=1}^\infty  k^m - 2 \sum_{k=1}^\infty  (2k)^m = (1-2^{1+m}) \sum_{k=1}^\infty k^m,
\label{connect} \end{equation}
relating the value of the  alternating series with the value of the non-alternating one.
But what of it? Surely, nothing would be gained by this move, as both series diverge.
Yes and no. 

Each alternating series still diverges in the usual sense, since even the necessary condition for convergence, that the general term $(-1)^{k+1} k^m$ must tend to zero as $k\to \infty$ is not met. However, as Euler knew full well \cite{euler1749}, these series {\em do\/} converge in
{\em the sense of Abel,\/} in today's parlance. This means that each generating function
$$ \sum_{k=1}^\infty  (-1)^{k+1} k^m x^k$$
has a finite limit as $x$ approaches $1$ from below along the real axis. The latter is not hard: indeed, if $m=0$, the resulting generating function is merely the geometric series $x/(1+x)=1-1/(1+x)$ and the other functions can be obtained from it by repeatedly using the differential operator $x{d \over dx}$: 
\begin{equation} \sum_{k=1}^\infty  (-1)^{k+1} k^m x^k = 1- \left( x {d\over dx} \right)^m {1\over 1+x}. \label{Abel} \end{equation}
From here, one can verify by induction that all alternating sums are indeed Abel-summable.
In view of this,  Euler made his `leap of faith': he {\em defined\/} each value $\zeta(-m)$  as the Abel sum of the alternating
 series $\sum_{k=1}^\infty (-1)^{k+1} k^m$ divided by the factor $(1-2^{1+m})$, as per~(\ref{connect}). 

Now, if we knew the Abel sums of the alternating series, we would know the values $\zeta(-m)$.
To get the first few, we can just apply~(\ref{Abel}) and plug in $x=1$
 (Euler lists more of these in \cite{euler1749}):
\begin{eqnarray*}
1^0 - 1^0 + 1^0  -\cdots & = & {1\over 2} \\
1^1-2^1+3^1 - \cdots & = & {1\over 4}  \\
1^2-2^2 +3^2 -\cdots & = & 0 \\
1^3-2^3 +3^3 -\cdots & = & {1\over 8}, \quad \text{etc.}
\end{eqnarray*}
So how can we evaluate all these Abel sums? In the spirit of Euler, let's make up an exponential generating function in a new symbol $z$ and manipulate it till we arrive at something useful:
\begin{eqnarray*} \sum_{m=0}^\infty {z^{m+1}\over m!} \left( \left( x {d\over dx} \right)^m {1\over 1+x}\right) \Big|_{x=1}  & = &
z \exp\left(zx {d\over dx} \right) \left( {1\over 1+x}\right) \Big|_{x=1}
\end{eqnarray*}
Note that the $m$th power of the operator $zx {d\over dx}$ maps each monomial $x^n$ to its multiple $ z^m n^m x^n$. Hence the exponential operator
$\exp\left(zx {d\over dx} \right)$ maps $x^n$ to $\sum_{m=0}^\infty {z^m n^m \over m!} x^n=\exp(zn) x^n$.  Hence, the same exponential operator maps the series
$1-x+x^2-\cdots=1/(1+x)$ to $\sum_{n=0}^\infty (-1)^n \exp (zn) x^n$. Finally, evaluating the last formula at $x=1$ yields the following generating function
$$ z \sum_{n=0}^\infty (-1)^n  \exp (zn) = {z\over 1+ e^z} ={z(e^z-1) \over e^{2z}-1}={z (e^z+1) -2z \over e^{2z} -1}={z\over e^z-1} -{2z\over e^{2z} -1}.$$
But this is just the difference between the original and the rescaled generating functions for the Bernoulli sequence!  So,
$$ {z\over e^z-1} -{2z\over e^{2z} -1} = \sum_{m=0}^\infty {B_m  (1-2^m)\over m!} z^m =\sum_{m=1}^\infty {B_m  (1-2^m)\over m!} z^m 
=\sum_{m=0}^\infty {B_{m+1}  (1-2^{m+1})\over (m+1)!} z^{m+1}  .$$
Comparing coefficients and recalling~(\ref{Abel}), we conclude
that the $m$th Abel sum $\sum_{k=1}^\infty (-1)^{k+1} k^m$ must be equal to the difference $-(1 - 2^{m+1})B_{m+1}/(m+1)$ for all $m\in \N$. In the special case $m=0$ the formula~(\ref{Abel}) yields the Abel sum $1+B_1=1/2=-B_1$. 

Since all odd Bernoulli numbers are zero except for $B_1$, our findings can be summarized as saying that the $m$th  Abel sum
$\sum_{k=1}^\infty (-1)^{k+1} k^m$ is equal to $(-1)^{m} (1-2^{m+1})B_{m+1}/(m+1)$. Hence, the value $\zeta(-m)$ should be taken to be
$(-1)^{m} B_{m+1}/(m+1)$, which is the familiar formula~(\ref{zetanegative}).

Apropos the odd Bernoulli numbers, the simplest way to see why they all vanish beginning with $B_3$ is to explicitly `peel away' the degree-one term from its generating function:
$$ {z\over e^z -1} - {B_1 z} = {z \over e^z - 1} +{z\over 2} ={z ( e^z +1) \over  2(e^z -1)} = {z ( e^{z/2} +e^{-z/2} ) \over  2(e^{z/2} - e^{-z/2})}.$$
The resulting function is now manifestly even, hence all its odd Taylor coefficients must be zero.

Back to the generating function used above to find the Abel sums of the alternating series, note that it is exactly the function $G(z)$ from Section~\ref{sec:euler} (up to a sign)! This shows that the values $\zeta(-m)$ encoded by $G(z)$ are precisely the same as the values 
derived from the alternating Abel sums. With our knowledge of today that the values $\zeta(-m)$ are well defined, Euler's  reduction of these values to the alternating Abel sums is no longer controversial. 

All this despite the fact that Euler himself actually evaluated $\zeta$ at nonpositive integers differently. To accomplish that, he combines the Euler-Maclaurin summation formula with the alternating trick to ascribe values to the  sums $x^m - (x+a)^m + (x+2a)^m -\cdots$. In other words, Euler replaces $a$ by $2a$ in~(\ref{EM}) to obtain the expansion of the even part $x^m + (x+2a)^m + (x+4a)^m + \cdots$,  which he then subtracts twice from $x^m + (x+a)^m + (x+2a)^m + \cdots$. This yields the expansion
$$ - \sum_{n=0}^\infty {a^{n-1}  B_n \over n!} f_m^{(n)}(x) + \sum_{n=0}^\infty { 2^n a^{n-1} B_n \over n!} f_m^{(n)}(x)  =\sum_{n=1}^\infty { (2^{n}-1) a^{n-1} B_n \over n!} f_m^{(n)}(x) $$ for the derivative of the `function' $x^m - (x+a)^m + (x+2a)^m -\cdots$ in terms of the derivatives of $f_m(x)=x^m$. Euler integrates this back to obtain the expansion 
$$ \sum_{n=0}^\infty { (2^{n+1}-1) a^{n} B_{n+1} \over (n+1)!} f_m^{(n)}(x) $$ for the alternating sum $x^m - (x+a)^m + (x+2a)^m -\cdots$
 and sets $a=1$ and 
$x=0$. Since $m$ is a nonnegative integer, all derivatives of $f$ evaluated at $0$ vanish unless the differentiation is performed exactly $m$ times. This yields  the value $${(2^{m+1}-1) B_{m+1} \over (m+1)!} \cdot m! = {(2^{m+1}-1) B_{m+1} \over (m+1)}$$ for the series $-1^m+2^m-3^m+\cdots$ and hence the value
$-B_{m+1}/(m+1)$ for the series $1^m+2^m+3^m+\cdots$, i.e., the value $\zeta(-m)$, $m\in \N$, again consistent with~(\ref{zetanegative}).

\section{Riemann $\leftrightarrow$ Euler translation}

I'd like to end this note with a very basic  Euler (= discrete) to Riemann (= continuous) translation.  

In the  discrete setting, the generating function $G$ must be differentiated $m$ times, then evaluated at $0$ to extract the value  of $\zeta$ at the nonpositive integer $-m$. In the continuous setting, the same task must be accomplished by integrating this same function against $z^{-m-1}$, creating a residue with the same value at the origin.  Thus the Hankel contour of integration, which separates the origin from the other singularities of the integrand, naturally arises. 

The `discrete' index $-m$ can now become a `continuous' complex variable $s$. The same Hankel contour instantly provides the analytic continuation for $\zeta(s)$ beyond the half-plane $\Re s >1$. The discrete functional relation naturally extends to the whole complex plane via the `inversion' of the contour and the residue theorem. One can point out that this relation extends to the complex plane simply because it  holds at all integer points, but the former strikes me as the `truer' reason. 

I hope this note will help those who wish to learn a bit about the Zeta function equipped only with undergraduate background in mathematics. In writing it, I was inspired by the monographs~\cite{borwein2008} and \cite{Hardy1963}, the lovely papers~\cite{Ayoub1974}, \cite{ball2022functional}, and -- especially  -- by the original works by Euler \cite{euler1731,euler1732/33, euler1734, euler1736,euler1737,euler1740,euler1749} and Riemann \cite{riemann1859uber}.  May the gentle reader be  encouraged, likewise, to read the classics by the words attributed to Laplace:
``{\em Read Euler, read Euler, he is the master of us all.}''

\section*{Acknowledgements}

I am grateful to Peter Sarnak, Keith Ball, and the anonymous referees for a number of valuable suggestions that helped improve exposition in this paper.

\small

\bibliographystyle{plainurl}
\bibliography{bib_zeta}


\end{document}